\documentclass{article}

\usepackage{mathrsfs,enumerate,amsmath}
\usepackage{amssymb}

 \newtheorem{thm}{Theorem}[section]
 \newtheorem{corr}[thm]{Corollary}

 \newtheorem{rem}[thm]{Remark}
 \newtheorem{conj}[thm]{Conjecture}
 \newtheorem{ex}[thm]{Example}
 \numberwithin{equation}{section}

\begin{document}

\title{R\'{e}nyi entropy and Tsallis entropy associated with positive linear operators}
\author{Ioan Ra\c{s}a\footnote{Department of Mathematics,
Technical University of Cluj-Napoca,
Memorandumului Street 28,
400114 Cluj-Napoca,
Romania, ioan.rasa$@$math.utcluj.ro}}
\date{}
\maketitle
%----------Author 1

Subjclass: 41A36, 94A17, 33C45


\begin{abstract}
This article is a continuation of my paper [arxiv: 1409.1015v2]. R\'{e}nyi and Tsallis entropies are associated to positive linear operators and properties of some functions related to these entropies are investigated.
\end{abstract}

\section{Introduction\label{sect:1}}

This paper is a continuation of~\cite{8}. In that article we considered discrete positive linear operators of the form
\begin{equation*}
Lf(x) = \sum _k f(x_k)a_k(x), \quad a_k(x)\geq 0, \quad \sum _k a_k(x)=1
\end{equation*}
for $x$ in some interval $I \subset \mathbb{R}$. In investigating the degree of non-multiplicativity of $L$ an important role was played by the function $S(x) = \sum _k a_k^2 (x)$; see~\cite{6}.

On the other hand, for each fixed $x\in I$ the numbers $(a_k(x))_k$ form a probability distribution. In this context $-\log S(x)$ is a R\'{e}nyi entropy~\cite{9} and $1-S(x)$ is a Tsallis entropy~\cite{10}. So the properties of the function $S(x)$ are relevant also in the study of these entropies.

Some properties of $S(x)$ were investigated in~\cite{5}-\cite{8}. In Section~\ref{sect:2} we continue to study such properties in the case of discrete operators. Section~\ref{sect:3} is devoted to some multivariate operators; see also~\cite{1}. In the last section we consider integral operators of the form $Lf(x)=\int _I K(x,t)f(t)dt$ with $K(x,t)\geq 0$, $\int _I K(x,t) dt = 1$, $x\in I$. Recall that in this case the associated R\'{e}nyi entropy is defined by $-\log \int _I K^2 (x,t)dt$, and the Tsallis entropy by $1-\int _I K^2 (x,t)dt$.

\section{A conjecture from~\cite{8} and some of its consequences\label{sect:2}}

We shall use the notation from~\cite{8}. In particular, we consider the function $S_{n,c}$ defined on the interval $I_c$. Conjecture 7.1 in~\cite{8} reads as follows:
\begin{equation*}
(C) \quad \log S_{n,c} \mbox{ is a convex function}.
\end{equation*}

Let us examine some consequences of this conjecture. With $X:=x(1+cx)$ and $X' = 1+2cx$, (3.10) in~\cite{8} can be written as
\begin{equation}
XX'S''_{n,c}(x) + (4(n+c)X+1)S'_{n,c}(x)+2nX'S_{n,c}(x)=0. \label{eq:2.1}
\end{equation}

$(C)$ is equivalent to $S''_{n,c}S_{n,c}\geq (S'_{n,c})^2$, and due to~\eqref{eq:2.1} both of them are equivalent to
\begin{equation*}
XX'\left ( \frac{S'_{n,c}}{S_{n,c}} \right )^2 + (4(n+c)X+1)\frac{S'_{n,c}}{S_{n,c}} + 2nX' \leq 0.
\end{equation*}

This leads immediately to
\begin{thm}\label{thm:2.1}
Conjecture $(C)$ is equivalent to
\begin{equation*}
(C') \quad \frac{S'_{n,c}(x)}{S_{n,c}(x)} \mbox{ is between $z_1(x)$ and $z_2(x)$},
\end{equation*}
where
\begin{equation*}
z_1(x) = \frac{-\sqrt{(1+4cX)^2+(4nX)^2}-(1+4cX)-4nX}{2XX'},
\end{equation*}
\begin{equation*}
z_2(x) = \frac{\sqrt{(1+4cX)^2+(4nX)^2}-(1+4cX)-4nX}{2XX'}.
\end{equation*}
\end{thm}

Let us consider the case $c=0$. Then the function $K_n(x):=S_{n,0}(x)$ is defined for $x\in [0,+\infty)$. In~\cite[(6.4)]{8} it was proved that
\begin{equation}
K_n(x)\leq \frac{1}{\sqrt{4nx+1}}, \quad x\geq 0. \label{eq:2.2}
\end{equation}

Recall also that $I_0$ is the modified Bessel function of first kind of order zero, and (see~\cite[(3.7)]{8})
\begin{equation}
I_0(x) = e^x K_n \left ( \frac{x}{2n} \right ), \quad x\geq 0. \label{eq:2.3}
\end{equation}

Combining~\eqref{eq:2.2} and~\eqref{eq:2.3} we get
\begin{equation}
I_0(x) \leq \frac{\exp{x}}{\sqrt{2x+1}}, \quad x\geq 0. \label{eq:2.4}
\end{equation}

\begin{corr}\label{corr:2.2}
Under the hypothesis that $\log S_{n,0}$ is convex, we have
\begin{equation}
\frac{-\sqrt{1+(4nt)^2}-1-4nt}{2t} \leq \frac{K'_n(t)}{K_n(t)} \leq \frac{\sqrt{1+(4nt)^2}-1-4nt}{2t}, \quad t> 0, \label{eq:2.5}
\end{equation}
\begin{equation}
K_n^2(x)\leq \frac{2\exp{\left ( \sqrt{1+(4nx)^2} -1-4nx \right )}}{\sqrt{1+(4nx)^2}+1}, \quad x\geq 0, \label{eq:2.6}
\end{equation}
\begin{equation}
I_0^2(x) \leq \frac{2\exp{\left ( \sqrt{1+4x^2} -1 \right )}}{\sqrt{1+4x^2} +1}, \quad x\geq 0. \label{eq:2.7}
\end{equation}
\end{corr}

\textbf{Proof.}
\eqref{eq:2.5} is a direct consequence of Theorem~\ref{thm:2.1}. \eqref{eq:2.6} can be obtained from the second inequality in~\eqref{eq:2.5} by integrating with respect to $t$ between $0$ and $x$. \eqref{eq:2.7} follows from~\eqref{eq:2.6} and~\eqref{eq:2.3}.
\newline \newline
In order to compare~\eqref{eq:2.2} with~\eqref{eq:2.6}, and~\eqref{eq:2.4} with~\eqref{eq:2.7}, it is easy to check that
\begin{equation*}
\frac{2\exp{\left ( \sqrt{1+(4nx)^2} -1-4nx \right )}}{\sqrt{1+(4nx)^2}+1} \leq \frac{1}{4nx+1}, \quad x\geq 0,
\end{equation*}

\begin{equation*}
\frac{2\exp{\left ( \sqrt{1+4x^2} -1 \right )}}{\sqrt{1+4x^2} +1} \leq \frac{\exp{2x}}{2x+1}, \quad x\geq 0.
\end{equation*}

Now let us consider the case $c=-1$. Then the function $F_n(x):=S_{n,-1}(x)$ is defined for $x\in [0,1]$. Let $x\in [0, 1/2)$ and $t = \frac{2x^2-2x+1}{1-2x}$. Then (see also~\cite[Section 4]{8}) $t\in [1,+\infty)$, $x = \frac{1-t+\sqrt{t^2-1}}{2}$, $X=x(1-x)=\sqrt{t^2-1}(t-\sqrt{t^2-1})/2$, $X'=1-2x=t-\sqrt{^2-1}$, $\frac{dt}{dx} = \frac{4X}{1-4X}$, $t=(1-2X)/X'$, $\sqrt{t^2-1}=2X/X'$.

Let $\left ( P_n(t) \right ) _{n\geq 0}$ be the Legendre polynomials. Then (see~\cite{new7}, \cite{7}, \cite[(4.6)]{8})
\begin{equation}
F_n(x) = \left ( t-\sqrt{t^2-1} \right )^n P_n(t). \label{eq:2.8}
\end{equation}

This yields
\begin{equation*}
\frac{F'_n(x)}{F_n(x)} = \left ( \frac{P'_n(t)}{P_n(t)} - \frac{n}{\sqrt{t^2-1}} \right ) \frac{dt}{dx},
\end{equation*}
and consequently
\begin{equation}
\frac{P'_n(t)}{P_n(t)} - \frac{n}{\sqrt{t^2-1}} = \frac{1-4X}{4X}\frac{F'_n(x)}{F_n(x)}. \label{eq:2.9}
\end{equation}

\begin{corr}\label{corr:2.3}
Under the hypothesis that $\log S_{n,-1}$ is convex, we have
\begin{equation}
\frac{F'_n(x)}{F_n(x)} \leq \frac{\sqrt{(1-4X)^2+(4nX)^2}-(1-4X)-4nX}{2XX'}, \quad x\in \left [ 0, \frac{1}{2} \right ), \label{eq:2.10}
\end{equation}
\begin{equation}
\frac{P'_n(t)}{P_n(t)} \leq \frac{\sqrt{4n^2(t^2-1)+(t-\sqrt{t^2-1})^2}-(t-\sqrt{t^2-1})}{2(t^2-1)}, \quad t>1. \label{eq:2.11}
\end{equation}
\end{corr}

\textbf{Proof.}
\eqref{eq:2.10} is a consequence of Theorem~\ref{thm:2.1}. \eqref{eq:2.11} follows from~\eqref{eq:2.10} and~\eqref{eq:2.9}.
\newline \newline
The following inequality was proved in~\cite[(1.2)]{7}:
\begin{equation}
\frac{n(n+1)}{2t+(n-1)\sqrt{t^2-1}} \leq \frac{P'_n(t)}{P_n(t)}, \quad t\geq 1. \label{eq:2.12}
\end{equation}

Using it and~\eqref{eq:2.9}, we get
\begin{equation}
-\frac{2nX'}{1+(n-3)X}\leq \frac{F'_n(x)}{F_n(x)}, \quad x\in \left [ 0, \frac{1}{2}\right ]. \label{eq:2.13}
\end{equation}

Other lower and upper bounds for $\frac{P'_n(t)}{P_n(t)}$ can be found in~\cite{7}. In particular, from~\cite[Theorems 2 and 3]{7} we have
\begin{equation}
\frac{P'_n(t)}{P_n(t)} \leq \frac{2n^2}{t+(2n-1)\sqrt{t^2-1}}, \quad t\geq 1, \label{eq:2.14}
\end{equation}
\begin{equation}
\frac{P'_n(t)}{P_n(t)} \leq \frac{n^2(2n+1)}{(n+1)t+(2n^2-1)\sqrt{t^2-1}}, \quad t \geq 1. \label{eq:2.15}
\end{equation}

\eqref{eq:2.11} and \eqref{eq:2.14} can be compared and we get
\begin{equation*}
\frac{\sqrt{4n^2 (t^2-1)+(t-\sqrt{t^2-1})^2}-(t-\sqrt{t^2-1})}{2(t^2-1)}\leq \frac{2n^2}{t+(2n-1)\sqrt{t^2-1}}, \quad t>1.
\end{equation*}

The inequality
\begin{equation*}
\frac{\sqrt{4n^2 (t^2-1)+(t-\sqrt{t^2-1})^2}-(t-\sqrt{t^2-1})}{2(t^2-1)} \leq \frac{n^2(2n+1)}{(n+1)t+(2n^2 - 1)\sqrt{t^2-1}}, \quad t>1,
\end{equation*}
is equivalent to
\begin{equation*}
\frac{t}{t+\sqrt{t^2-1}}\geq \frac{3n+2}{4n+3}.
\end{equation*}

This last inequality is true for $t$ approaching $1$, and false for $t$ approaching $+\infty$.

Let us remark that~\eqref{eq:2.15} yields by integration
\begin{equation}
P_n(t)\leq (t+\sqrt{t^2-1})^{\frac{n(2n^2-1)}{2n^2-n-2}}\left ( t + \frac{2n^2-1}{n+1} \sqrt{t^2-1} \right ) ^{-\frac{n(n+1)}{2n^2-n-2}}, \quad t\geq 1. \label{eq:2.16}
\end{equation}

This inequality is stronger than
\begin{equation}
P_n(t) \leq (t+\sqrt{t^2-1})^{\frac{n(2n-1)}{2(n-1)}} \left(t+(2n-1)\sqrt{t^2-1}\right)^{-\frac{n}{2(n-1)}}, \quad t\geq 1, n\geq 2, \label{eq:2.17}
\end{equation}
which can be obtained from~\eqref{eq:2.14}.

We conclude this section with a remark concerning the function $A_{n,c}:=\frac{S'_{n,c}}{S_{n,c}}$. By using~\cite[(3.10)]{8}, or~\eqref{eq:2.1}, we deduce easily that $A_{n,c}$ satisfies the Riccati equation
\begin{equation*}
x(1+cx)(1+2cx) (A'_{n,c}+A^2_{n,c}) + (4(n+c)x(1+cx)+1)A_{n,c}+2n(1+2cx)=0.
\end{equation*}

\section{Multivariate operators\label{sect:3}}

First, consider the classical Bernstein operators on the canonical simplex of $\mathbb{R}^2$. The sum of the squared fundamental Bernstein polynomials is in this case
\begin{equation*}
R_n(x,y):= \sum_{i+j\leq n} \left ( \frac{n!}{i!j!(n-i-j)!} \right ) ^2 x^{2i}y^{2j}(1-x-y)^{2(n-i-j)}
\end{equation*}
\begin{equation*}
=\sum _{j=0}^n \sum _{i=0}^{n-j} {n \choose j}^2 {n-j \choose i}^2 x^{2i} y^{2j}(1-x-y)^{2(n-i-j)},
\end{equation*}
for $x\geq 0$, $y\geq 0$, $x+y \leq 1$; see~\cite[(6.3.6)]{2}, \cite[Sect. 3.1.2]{3}.

Let $y\in [0,1)$ be fixed. Then for $x\in [0,1-y]$ we have
\begin{equation*}
R_n(x,y) = \sum _{j=0}^n {n \choose j}^2 y^{2j}(1-y)^{2(n-j)}\sum _{i=0}^{n-j} {n-j \choose i}^2 \left ( \frac{x}{1-y} \right )^{2i} \left ( 1-\frac{x}{1-y} \right ) ^{2(n-j-i)}.
\end{equation*}

For each $j\in \{ 0,1, ..., n \}$,
\begin{equation*}
\sum _{i=0}^{n-j} {n-j \choose i}^2  \left( \frac{x}{1-y} \right )^{2i}\left ( 1-\frac{x}{1-y} \right ) ^{2(n-j-i)} = F_{n-j} \left( \frac{x}{1-y} \right ),
\end{equation*}
where $F_{n-j}=S_{n-j,-1}$. It is known (see~\cite{5}, \cite{6}, \cite{7}, \cite{8}) that $F_{n-j}$ is convex on $[0,1]$. It follows that for each fixed $y\in [0,1)$, the function $R_n(\cdot , y)$ is convex on $[0, 1-y]$. In other words, $R_n$ is convex on each segment parallel to $Ox$. Similarly we see that $R_n$ is convex on each segment parallel to a side of the canonical triangle of $\mathbb{R}^2$. This means that $R_n$ is \emph{axially-convex}; concerning this terminology see~\cite[p. 407]{2}, \cite[Sect. 3.5]{3}.

Now consider the classical Bernstein operators on the square $[0,1]^2$: see~\cite[(6.3.101)]{2}, \cite[Sect. 3.1.5]{3}. The sum of the squared fundamental Bernstein polynomials is in this case
\begin{equation*}
Q_n(x,y) = \sum _{i=0}^n \sum _{j=0}^n {n \choose i}^2 {n \choose j}^2 x^{2i} (1-x)^{2n-2i} y^{2j} (1-y)^{2n-2j} = F_n(x)F_n(y).
\end{equation*}

It is easy to verify that the following three statements are equivalent:
\begin{enumerate}[i)]
\item{} $\log F_n$ is convex on $[0,1]$;
\item{} $Q_n$ is convex on $[0,1]^2$;
\item{} $\log Q_n$ is convex on $[0,1]^2$.
\end{enumerate}

\section{Entropy and variance. Integral operators}\label{sect:4}

Let $I$ be an interval and $L$ a positive linear operator on a space of functions defined on $I$, containing the functions $e_i(x)=x^i$, $i=0,1,2$. Suppose that $Le_0=e_0$.

The \emph{variance} associated with $L$ is the function
\begin{equation*}
V(x):=Le_2(x)-(Le_1(x))^2, \quad x\in I.
\end{equation*}

If $L$ is a discrete operator of the form $Lf(x) = \sum _k f(x_k)a_k(x)$, let $S(x):=\sum _k a_k^2(x)$, $x\in I$. If $L$ is an integral operator of the form $Lf(x) = \int _I K(x,t)f(t)dt$, let $S(x):= \int _I K^2(x,t)dt$, $x\in I$.

In both cases the R\'{e}nyi entropy associated with $L$ is $-\log S(x)$, and the Tsallis entropy is $1-S(x)$, $x\in I$.

\begin{ex}\label{ex:4.1}
\end{ex}
Let $L_{n,c}f(x) = \sum _{j=0}^\infty f \left ( \frac{j}{n} \right ) p_{n,j}^{[c]} (x)$, see \cite[Sect. 2]{8}. Then $S(x) = S_{n,c}(x)$ and $V(x) = V_{n,c}(x) = \frac{x(1+cx)}{n}$. According to \cite[(3.5), (3.8)]{8},
\begin{equation}
S_{n,c}(x) = \frac{1}{\pi} \int _0 ^\pi \left ( 1+4ncV_{n,c}(x)\sin ^2 \frac{\varphi}{2} \right )^{-n/c}d\varphi, \quad c \neq 0, \label{eq:4.1}
\end{equation}
\begin{equation}
S_{n,0}(x) = \frac{1}{\pi} \int _0 ^\pi \exp{\left ( -4n^2 V_{n,0}(x)\sin ^2 \frac{\varphi}{2} \right )}d\varphi .\label{eq:4.2}
\end{equation}

\begin{ex}\label{ex:4.2}
\end{ex}
For the Kantorovich operators~\cite[p. 333]{2} we have $S_n(x) = (n+1)F_n(x)$ and $V_n(x)=(n+1)^{-2} \left ( nx(1-x) + \frac{1}{12} \right )$.

\begin{ex}\label{ex:4.3}
\end{ex}
Consider the Gauss-Weierstrass operators~\cite[p. 310]{2}, \cite[p. 114]{4}:
\begin{equation*}
W_r f (x) = \int _\mathbb{R} (4\pi r)^{-1/2} \exp{\left ( - \frac{(t-x)^2}{4r} \right )} f(t)dt, \quad r>0.
\end{equation*}

Then $V_r(x)=2r$ and $S_r(x) = (8\pi r)^{-1/2}$, $x\in \mathbb{R}$.

Generally speaking, for a convolution operator
\begin{equation*}
Lf(x)=\int _\mathbb{R} \varphi (x-t)f (t)dt \mbox{ we have } V(x) = \int _\mathbb{R} s^2 \varphi (s)ds - \left ( \int _\mathbb{R} s \varphi (s)ds \right )^2
\end{equation*}
and $S(x) = \int _\mathbb{R} \varphi ^2(s)ds$, so that $V$ and $S$ are constant functions.

\begin{ex}\label{ex:4.4}
\end{ex}
For the Post-Widder operators~\cite[p. 114]{4},
\begin{equation*}
V_n(x) = \frac{x^2}{n} \mbox{ and } S_n(x) = {2n-2 \choose n-1} 2^{1-2n}\frac{n}{x}, \quad x>0.
\end{equation*}

\begin{ex}\label{ex:4.5}
\end{ex}
Consider the Durrmeyer operators~\cite[p. 335]{2}.

In this case
\begin{equation*}
V_n(x) = \frac{n+1}{(n+2)^2(n+3)} (2nx(1-x)+1),
\end{equation*}
\begin{equation*}
S_n(x) = \sum _{k=0}^{2n} c_{n,k} {2n \choose k}x^k (1-x)^{2n-k},
\end{equation*}
where
\begin{equation*}
c_{n,k}:=\frac{(n+1)^2}{2n+1} {2n \choose k}^{-2}\sum _{j=0}^k {n \choose j}^2 {n \choose k-j}^2, \quad k = 0, 1, ... , 2n,
\end{equation*}
where, as usual, ${n \choose m} = 0$ if $m>n$.

It is easy to see that $c_{n,2n-k} = c_{n,k}$, $k=0,1,..., 2n$.

\begin{conj}\label{conj:4.6}
The sequence $(c_{n,k})_{k=0,1,..,2n}$ is convex and, consequently, the function $S_n$ is convex on $[0,1]$.
\end{conj}

\begin{ex}\label{ex:4.7}
\end{ex}
For the genuine Bernstein-Durrmeyer operators, defined by
\begin{equation*}
U_nf(x) = f(0)b_{n,0}(x) + f(1)b_{n,n}(x)+(n-1)\sum _{k=1}^{n-1}b_{n,k}(x)\int _0^1 b_{n-2,k-1}(t)f(t)dt,
\end{equation*}
with $b_{n,k} (x)= {n \choose k}x^k (1-x)^{n-k}$, we have
\begin{equation*}
V_n(x) = \frac{2x(1-x)}{n+1}
\end{equation*}
and
\begin{eqnarray*}
&& S_n(x) = (1-x)^{2n}+x^{2n} + \\ &&+ \frac{(n-1)^2}{2n-3} \sum _{k,j=1}^{n-1} {n-2 \choose k-1} {n-2 \choose j-1}{n \choose k}{n \choose j}{2n-4 \choose k+j-2}^{-1} x^{k+j}(1-x)^{2n-k-j}.
\end{eqnarray*}

\begin{rem}
In Examples~\ref{ex:4.1}-\ref{ex:4.4}, and also in Example~\ref{ex:4.5} under Conjecture~\ref{conj:4.6}, the functions $V(x)$, $1-S(x)$ and $-\log S(x)$ are all increasing or all decreasing on suitable subintervals of $I$. In other words, the variance, the Tsallis entropy and the R\'{e}nyi entropy are synchronous functions.
\end{rem}

\subsubsection*{Acknowledgement}

The author is grateful to Dr. Gabriela Raluca Mocanu for inspiring discussions.
% ------------------------------------------------------------------------

% ------------------------------------------------------------------------
\end{document}